\documentclass[preprint,12pt]{elsarticle}

\usepackage{units}

\usepackage{amssymb}
\usepackage{amsmath}

\newcounter{bla}

\journal{journal for publication}

\usepackage{graphicx}
\usepackage[normalem]{ulem}
\usepackage{color}
\usepackage{dcolumn}
\usepackage{bm}
\usepackage{float}

\begin{document}

\begin{frontmatter}



\title{Transforming Stiffness and Chaos}


\author{Jan Scheffel \corref{author}}

\cortext[author] {Jan Scheffel.\\\textit{E-mail address:} jans@kth.se}
\address{Electromagnetic Engineering and Fusion Science \\ 
	KTH Royal Institute of Technology, Stockholm, Sweden}

\begin{abstract}
Stiff and chaotic differential equations are challenging for time-stepping numerical methods. For explicit methods, the required time step resolution significantly exceeds the resolution associated with the smoothness of the exact solution for specified accuracy. In order to improve efficiency, the question arises whether transformation to asymptotically stable solutions can be performed, for which neighbouring solutions converge towards each other at a controlled rate. Employing the concept of local Lyapunov exponents, it is demonstrated that chaotic differential equations can be successfully transformed to obtain high accuracy, whereas stiff equations cannot. For instance, the accuracy of explicit fourth order Runge-Kutta solution of the Lorenz chaotic equations can be increased by two orders of magnitude. Alternatively, the time step can be significantly extended with retained accuracy. \end{abstract}


\end{frontmatter}

\section{Introduction}

\noindent When assessing the usefulness of a computational method for differential equations, there are at least five properties of the problems addressed that come into play: 1) \textit{stability} (characterized by Lyapunov exponents $\lambda_i < 0$, with small magnitudes), 2) \textit{oscillations} (number of extrema in the computational interval), 3) \textit{steepness} (normalized maximum rate of change), 4) \textit{stiffness} (at least one $\lambda_i < 0$, with large magnitude), and 5) \textit{chaoticity} (continuous, non-periodic solution in at least three variables \cite{Hirsch,Sprott}, with at least one $\lambda_i > 0$). We will, in this work, focus on the concepts of stiffness and chaoticity for explicit time-stepping methods, making frequent use of the concept of local Lyapunov exponents.

The phenomenon of stiffness may severly limit the time steps of explicit methods. The concept has been found difficult to characterize precisely but substantial progress was made in the work by Cartwright \cite{Cartwright:1}.  It is generally acknowledged that for stability a stiff differential equation requires step lengths that are "excessively small in relation to the smoothness of the exact solution" in the interval \cite{Lambert}. The reason for this behaviour is that the scales for decay of different parts of the solution are widely separated. To handle stiffness within time-stepping methods one must usually resort to time-consuming, less accurate implicit or semi-implicit methods. An alternative is time-spectral methods \cite{Scheffel:GWRM2}. 

The notion of Lyapunov exponents is helpful when distinguishing between stiffness and chaoticity for nonlinear problems. A novel method for approximate computation of local Lyapunov exponents (LLEs), without the necessity of solving the differential equations, is presented in Section 2. By "local" is here meant that the exponents are associated with the current time step or time interval. We also present an improved approach for categorization and formal separation of the properties of stiffness and chaoticity, in line with that of \cite{Cartwright:1}. LLEs are applied to examples of stiff linear and nonlinear equations in Sections 3 and 4, illustrating their usefulness. 

In Sections 5 and 6, the concept of LLEs is instrumental for proving that chaotic differential equations can be transformed into differential equations with asymptotically stable solutions, and solved efficiently at strongly improved accuracy. It is also shown that stiff differential equations cannot be transformed similarly since the difficulty to numerically resolve neighbouring, rapidly diminishing solutions, is transferred to the numerical resolution of the transformed functions. Notably, time-spectral methods \cite{Scheffel:GWRM2} remain unaffected by stiffness or chaos, as they are methodologically acausal, or in a sense, implicit. Consequently, they often represent optimal choices for addressing stiff or chaotic problems, a topic that will be elaborated upon in a forthcoming publication.

A discussion and summary of the transformation method proposed in this work can be found in Sections 7 and 8.

\section{Determining local stiffness and chaoticity}
\noindent It is of interest to find a formal and practical characterization for determining stiffness and chaoticity of a given linear or nonlinear system of differential equations, both locally and globally. The concept of Lyapunov exponents is useful for separating the properties of stability, stiffness and chaoticity. We here employ a standard definition for global Lyapunov exponents, whereafter we move on to present some novel aspects for LLEs. 

Assume that a coupled set of $N$ linear or nonlinear differential equations are solved for two different initial conditions with initial, infinitesimal separation $\bf{\delta {X}_0}$ (components $\delta X_{0i}$) in phase space. The time-dependent separation of solution curves can be approximated within a linear analysis;

\begin{equation}
	|\delta X_{i}(t)| \approx e^{\lambda_i t} |\delta X_{0i}|  .
\end{equation}

\noindent This relation defines the spectrum of $N$ Lyapunov exponents $\lambda_i$, determined by separations in $N$-dimensional space. The largest of these Lyapunov exponents is termed \textit{the maximal Lyapunov exponent} $\lambda_M$. For a \textit{chaotic} problem, $\lambda_M$ is positive, thus initially closely separated orbits will diverge exponentially in time. For many problems, however, $\lambda_M$ is instead moderately negative so that the solution is \textit{stable} in the sense that perturbed solutions remain close to the original solution, or even \textit{asymptotically stable} in which case original and perturbed solutions converge towards each other over time \cite{Heath}. The latter property saves convergence and accuracy even for basic, low order explicit difference methods with large time steps.

For \textit{stiff} problems, the minimum Lyapunov exponents are strongly negative. The associated rapid convergence may be difficult to resolve for finite difference methods. Whether a problem is stiff or not requires an analysis of the local variation of the solution, since stiffness implies that the local step size must be much smaller than what is necessitated by numerical resolution of the solution. For scalar functions $u=u(t)$, the absolute local curvature $\kappa(t)$ is a suitable measure of the local variation;

\begin{equation}
	\kappa(t) =  \frac{|u''(t)|}{(1+{u'(t)}^2)^{3/2}} \quad .
\end{equation}

\noindent \noindent A prime denotes differentiation with respect to $t$. \textit{Local Lyapunov exponents} (LLEs)  \cite{Cartwright:1} can now be defined for each phase space dimension $i$;

\begin{equation}
	\gamma_i(\tau;t)= \lim_{\sigma_i(\tau)\rightarrow 0} \frac{1}{\tau} ln \frac{\sigma_i(t+\tau)}{\sigma_i(t)}  \quad .
\end{equation}

\noindent \noindent Here $\sigma_i(\tau)$ denotes the $N$ principal axes of a small ellipsoid in phase space, representing the temporally varying separations. By comparing with the standard definition of global Lyapunov exponents (compare Eq. (1))

\begin{equation}
	\lambda_i= \lim_{t \rightarrow \infty} \,  \lim_{|\delta X_{0i}| \rightarrow 0} \, \frac{1}{t} ln \frac{|\delta X_i(t)|}{|\delta X_{0i}|}  \quad ,
\end{equation}

\noindent \noindent it is clear that local Lyapunov exponents, for any time $t$, estimate the deviation of the original and perturbed solutions in an interval $[t,t+\tau]$, or equivalently $t \in [t_n,t_{n+1}]$ with $t_{n+1}-t_n=\tau$. Lyapunov exponents are obtained from local Lyapunov exponents from the relation

\begin{equation}
	\lambda_i = \lim_{\tau \rightarrow \infty} \gamma_i(\tau;t)  \quad .
\end{equation}

\noindent \noindent Following Cartwright \cite{Cartwright:1}, stiffness can be defined for nonlinear systems, still allowing for chaoticity:  \textit {A system is stiff in a given interval if in that interval the most negative local Lyapunov exponent is large, while the curvature of the trajectory is small.} This criterion can be employed for stiff equations using the measure
\begin{equation}
	R = \frac{|\gamma_{i_{min}}(\tau;t)|}{\kappa(t)}
\end{equation}

\noindent \noindent where $\gamma_{i_{min}}$ denotes the most negative of the $n$ local Lyapunov exponents. Since at least one positive Lyapunov exponent is required for chaos and at least one should be strongly negative for stiffness, \textit{stiff chaos} implies a large spread of Lyapunov exponents. The measure is preferable to more pragmatic earlier definitions like "Stiff equations are problems for which explicit methods don't work" \cite{Hairer}. It also extends the linear concept of stiffness ratio \cite{Lambert} to nonlinear problems. Furthermore, the definition avoids the assumption that stiff problems should be non-chaotic, as in Lambert's definition \cite{Lambert}.

The definition enables distinction between, for example, systems that are stiff but non-chaotic ($R$ large and $\gamma_{i_{max}} \leq0$) and systems that are both stiff and chaotic ($R$ large and $\gamma_{i_{max}} > 0$). If the curvature associated with stiffness  $\kappa_{stiff} (t)$ is comparable to, or exceeds, the curvature $\kappa (t)$ of $u(t)$, accuracy is determined by stiffness rather than by the functional variation of $u(t)$. 

The maximum step length ${\Delta t}_{max}$ for approximating a function $u(t)$ with local curvature $\kappa(t)$ at accuracy $\epsilon$  in $t \in [t_n,t_{n+1}]$ can be estimated as follows. A small segment of a local curvature radius $1/\kappa(t)$ can be approximated to accuracy $\epsilon$ by a secant with length $\Delta t$ as a function of $t$. A simple geometric analysis straightforwardly yields the approximate relation 

\begin{equation}
	(\Delta t)_{max} = 2 \sqrt{2} \sqrt{\frac{\epsilon}{\kappa(t^*)}}
\end{equation} 

\noindent Note that $(\Delta t)_{max}$ is a function of time $t^*=t-t_n$. 

Next we determine the maximum step size $(\Delta t)_{stiff} $ that is allowed by stiffness in $t \in [t_n,t_{n+1}]$, for specified accuracy $\epsilon$. The focus is on contributions to the solution, or neighbouring solutions, in this interval, with local dependence of the form $\delta u (t) = \epsilon e^{\gamma t^*}$ where $\gamma$ is negative for stiff problems. Since $\epsilon$ is the desired accuracy, which should be maintained throughout the solution algorithm, the choice of $\delta u(0)$ is reasonable. In the general case, $\gamma=\gamma(t)$, but we will here assume that $\gamma$ is a constant, that may however differ for each interval $[t_n,t_{n+1}]$. If the differential equation is stiff, terms like $\epsilon e^{\gamma t^*}$ with associated curvature 

\begin{equation}
	\kappa_{stiff}(t^*) =  \frac{\epsilon \gamma^2 e^{\gamma t^*}}{(1+\epsilon^2 \gamma^2 e^{2 \gamma t^*})^{3/2}}  
\end{equation}

\noindent must be resolved; it does not suffice to resolve $u(t)$ itself to maintain $\epsilon$ accuracy overall. From Eq. (8), it is seen that the maximum value of $\kappa_{stiff}$ is independent of $\epsilon$; $\kappa_{stiff,max}=(2\sqrt{3} /9) |\gamma|$ and occurs for $t^*_{max}=-ln(2 \gamma^2 \epsilon^2)/(2 \gamma)$.

At this time (for $t^*_{max}>0$) the maximum allowable time step for accuracy $\epsilon$ can be obtained using Eq. (7);

\begin{equation}
	(\Delta t)_{stiff} = 6 \sqrt{\frac{\epsilon}{\sqrt{3} |\gamma|}}
\end{equation} 

\noindent In most cases $2 \gamma^2 \epsilon^2 < 1$ however, and $t^*_{max} < 0$, so that the maximum value of $\kappa$ is obtained at the start $t^*=0$ of the interval $[t_n,t_{n+1}]$;

\begin{equation}
	\kappa_{stiff,max}=\frac{\epsilon\gamma^2}{(1+\epsilon^2 \gamma^2)^{3/2}} .
\end{equation} 

\noindent For this case the maximum allowed step length is found at $t^*=0$ and is

\begin{equation}
	(\Delta t)_{stiff} = \frac{2 \sqrt{2}}{|\gamma|}(1+\epsilon^2 \gamma^2)^{3/4}   .
\end{equation} 

Thus we have estimated the maximum allowed step length $(\Delta t)_{stiff}$, determined by the stiffness of the problem, to be compared with $(\Delta t)_{max}$, required by the specific numerical method for resolving the solution curve to desired accuracy. For $(\Delta t)_{max} \leq (\Delta t)_{stiff}$, the differential equation should not be regarded as stiff. In adaptive methods, $(\Delta t)_{max}$ and $(\Delta t)_{stiff}$ can be continually evaluated for determining local stiffness or chaoticity. When, for example, applying the transformation techniques described in Section 5, a more intuitive measure of stiffness than the ratio $R$ of Eq. (6) is useful:

\begin{equation}
	Q(t) \equiv \frac{(\Delta t)_{max}}{(\Delta t)_{stiff}}  > 1,
\end{equation} 

\noindent A solution algorithm for which $Q(t)<1$ is thus not limited by stiffness. Note that ${(\Delta t)}_{max}$ is given by Eq. (7) and ${(\Delta t)}_{stiff}$ is found from either of Eqs. (9), for $2 \gamma^2 \epsilon^2 \geq 1$, or Eq. (11), for $2 \gamma^2 \epsilon^2 < 1$. This formalism will now be applied to stiff and chaotic problems.

\section{A stiff linear problem}
\noindent The success of numerical integration techniques relies to a large extent on that neighbouring solutions are asymptotically stable, converging towards each other over time and damping out numerical errors. Stiffness appears, however, when this convergence becomes extremely rapid, producing strong local gradients that upset explicit schemes unless very small time steps are used. We now solve a typical, stiff linear problem \cite{Heath} and use the formalism of Section 2 to diagnose its stiffness properties;

\begin{equation}
	\frac{du}{dt} = -au+at+a+1, \quad u(0) = 1.
\end{equation}

The general solution is $u(t)=1+t+ce^{-at}$, with $c$ being an arbitrary constant. The solution satisfying the given initial condition is $u(t)=1+t$. For $a \gg 1$ and large $t$ the general solution becomes a straight line, independent of the initial condition. For $u(0)=1$ the solution is smooth (since $c=0$), but for neighbouring initial conditions (other values of $c$) the solution changes very rapidly for small $t$, being problematic for explicit methods. The cause of this behaviour is the strongly negative Lyapunov exponent, due to the large value of $a$, in combination with the extremely small trajectory curvature (\cite{Cartwright:1,Heath}). For this problem, we have access to the exact solution, which shows that even small numerical errors will make the equation stiff. 

\subsection{Approximate method for estimating local stiffness}

\noindent We now introduce an approximative formal method to determine the character, being it stable, stiff or chaotic, of linear and non-linear initial-value problems, using LLEs. Other methods for determining LLEs are proposed in the literature \cite{Ayers}, but often feature higher computational cost due to the need for computing long system trajectories.

Denoting the exact solution with $u(t;0)$ and an $\epsilon$-perturbed ($\epsilon \ll 1$) solution $u(t;\epsilon)$, obtained by solving Eq. (13) for $u(0)=1+\epsilon$ instead, the deviation of the solutions becomes
\begin{equation}
	\delta u(t) = u(t; \epsilon) - u(t;0) .
\end{equation}

\noindent Subtracting the two versions of Eq. (13) for $u(t;\epsilon)$ and $u(t;0)$, there results
\begin{equation}
	\frac{d \delta u}{dt} = -a \delta u
\end{equation}

\noindent The solution to Eq. (15), for each time interval, is
\begin{equation}
	\delta u(t)=\delta u (t_n)e^{-at^*}
\end{equation}

\noindent where $t \in [t_n,t_{n+1}]$, with $t_n$ arbitrary and $t^*=t-t_n$. Hence the solution to Eq. (13) is asymptotically stable for moderately large $a>0$ and stiff for $a \gg 1$. To study stiffness, we perturb the initial condition $u(0;\epsilon) = 1+ \epsilon = 1.05$ and set $a=300$. The perturbed solution $u(t;\epsilon)=1+t+0.05e^{-300t}$ deviates marginally, but rapidly from $u(t;0)=1+t$ at low $t$. 

\begin{figure}[h!]
	\centering
	\includegraphics[width=5in]{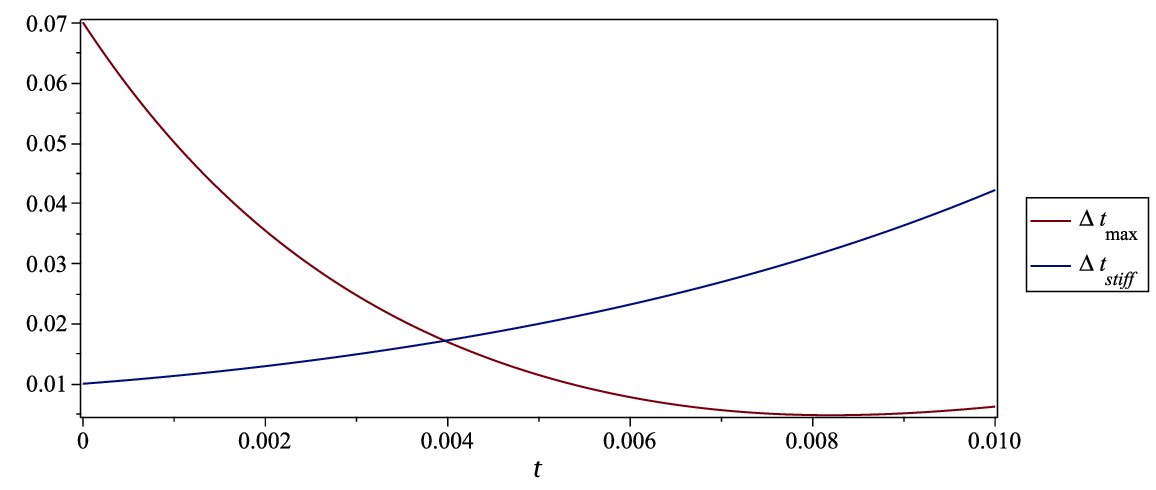}
	\caption{Comparison, for solution of Eq. (13) at $\epsilon=0.001$ accuracy, between maximum time step $(\Delta t)_{max}$ for numerical resolution and maximum time step $(\Delta t)_{stiff}$ for resolving stiffness of perturbed solution $u(t;\epsilon)=1+t+0.05e^{-300t}$. It is seen that the ODE is locally stiff for $t<0.004$, where $Q>1$. }
\end{figure}

The ODE is not stiff everywhere; see Figure 1. For $\delta u(t_0) = \epsilon = 0.001$ there results $2\epsilon^2 a^2=0.18 < 1$, indicating that the maximal stiffness obtains at the beginning of the  interval, and $\kappa_{stiff,max} \approx 90$. This renders $(\Delta t)_{stiff}=0.0094$. Using Eqs. (7) and (11), $(\Delta t)_{max}$ is computed as function of $t$. The problem is not stiff when the curvature of $u(t;\epsilon)$ is sufficiently large, so that $Q(t)=(\Delta t)_{max} /(\Delta t)_{stiff}<1$, which occurs for $t>0.004$, as can be seen in Figure 1. Apparently the condition $Q(t)>1$ for stiffness is an intuitively transparent and straightforward diagnostic.

Somewhat surprisingly, fourth order Runge-Kutta (RK4) solves this stiff problem well. Stability can be reached already using as few as 12 time steps, although at a relative error of $0.03$. An accuracy of $0.001$ is obtained using 25 time steps. Solution of Eq. (13) using an implicit time-stepping method (second order Trapezoidal method) \cite{Scheffel:GWRM3}, requires some 40 time steps for this accuracy. 

In summary we have, using a simple model problem, seen how LLEs can be used to define and compute the time-dependent diagnostics $\Delta t_{max}$, $\Delta t_{stiff}$ and $Q(t)$ for assessing local stiffness of an ODE.

\section{Stiff nonlinear problems}
\noindent We next consider \textit{nonlinear} stiff problems. Two example problems will be studied; a differential equation for modelling explosive combustion and a problem in kinetic analysis of chemical reactions. 

\subsection{Explosive combustion}
\noindent The dynamic evolution of a flame can, in spherically symmetric geometry, be modelled \cite{Scheffel:GWRM1} with the equation 

\begin{equation}
	\frac{du}{dt} = u^2-u^3, \quad u(0) = d.
\end{equation}

Here $d\in ]0,1[$ is an arbitrary constant. This ODE is stiff for $t>1/d$ when $d$ is small. Solutions for $d=0.1$, $d=0.05$ and $d=0.01$ are shown in Fig. 3.

\begin{figure}[h!]
	\centering
	\includegraphics[width=5.5in]{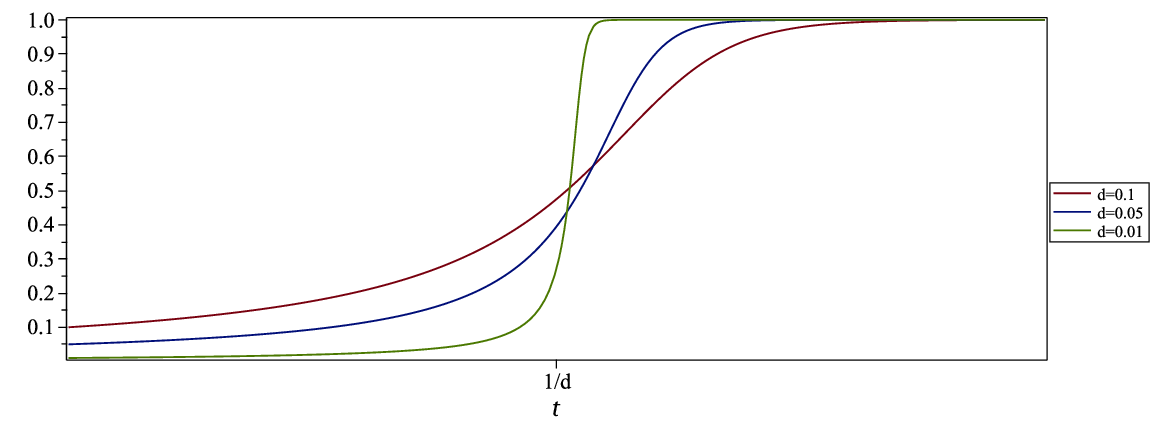}
	\caption{GWRM solutions $u(t)$ of Eq. (17) for $d=0.1$, $0.05$ and $0.01$.}
\end{figure}

How can the claim that the equation is stiff be justified? We employ again the method described in the previous section. Assuming the existence of an $\epsilon$-perturbed solution to Eq. (17), we label the two solutions $u(t;0)$ and $u(t;\epsilon)$, letting $\delta u(t)=u(t;\epsilon)-u(t;0)$, and obtain 

\begin{equation}
	\frac {d \delta u (t)}{dt} = \frac {du (t;\epsilon)}{dt}-\frac {du (t;0)}{dt}=-\delta u f(u(t;0))
\end{equation}

\noindent to first order in $\delta u$, where $f$ is a function to be found. By subtraction of the two versions of Eq. (17) we indeed have

\begin{equation}
	\frac{d \delta u}{dt} = -\delta u (3u^2(t;0)-2u(t;0)) .
\end{equation}

\noindent As $u(t;0) \rightarrow 1$, the solution to Eq. (19) becomes $\delta u = \delta u_{0}e^{-(t-t_0)} = \delta u_{0}e^{-\tau}$, which explains stiffness when $t \rightarrow 2/d >> 1$ for small $d$ in this interval. If a constant step length $\Delta t$ is employed throughout, it will for $d \leq 0.05$  be decided by $(\Delta t)_{max}$ (obtained from Eqs. (7) and (17)) near $t=1/d$ and by stiffness through $(\Delta t)_{stiff}$ (from Eqs. (11) and (19)) for $t>1/d$. 

\begin{figure}[h!]
	\centering
	\includegraphics[width=2.5in]{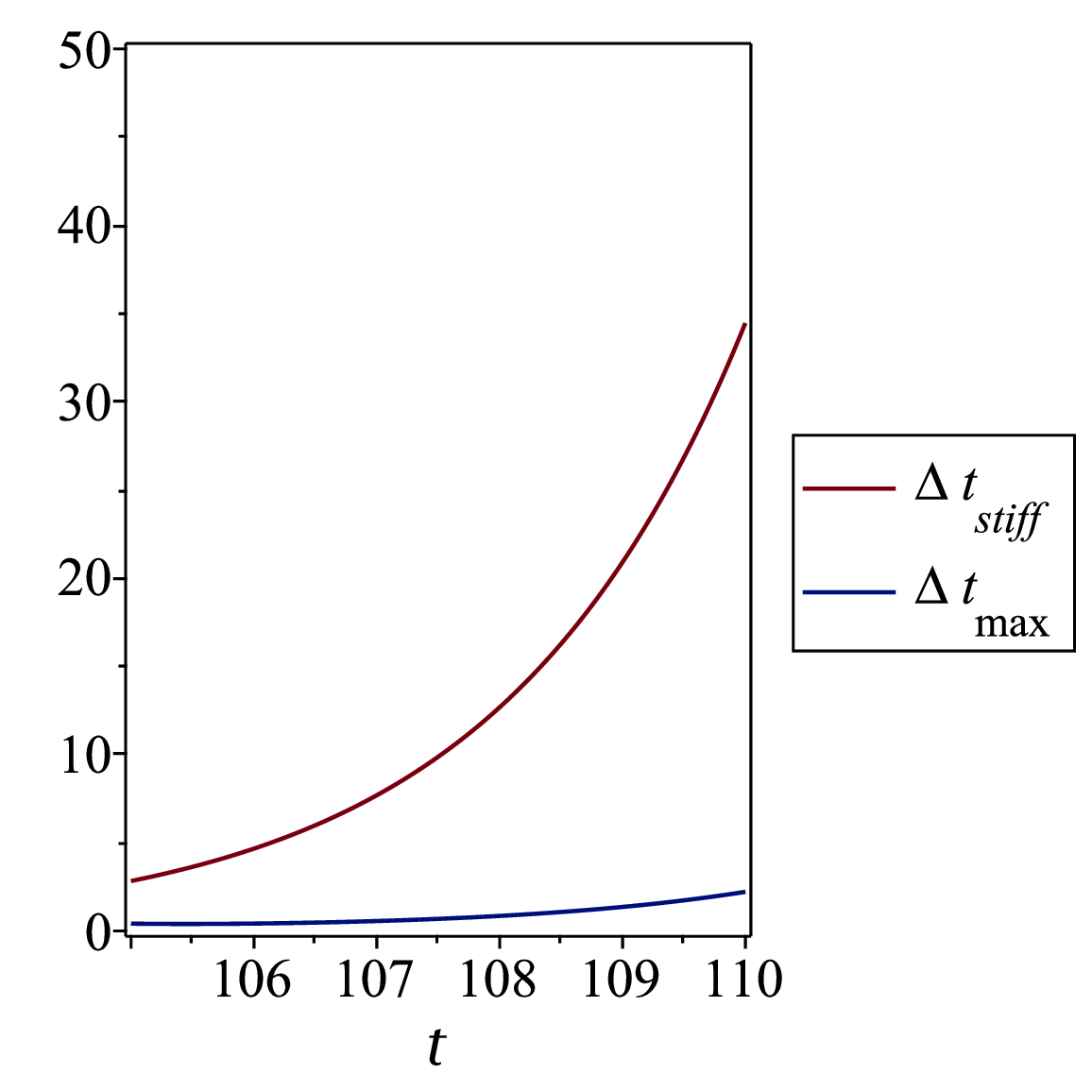}
	\includegraphics[width=2.5in]{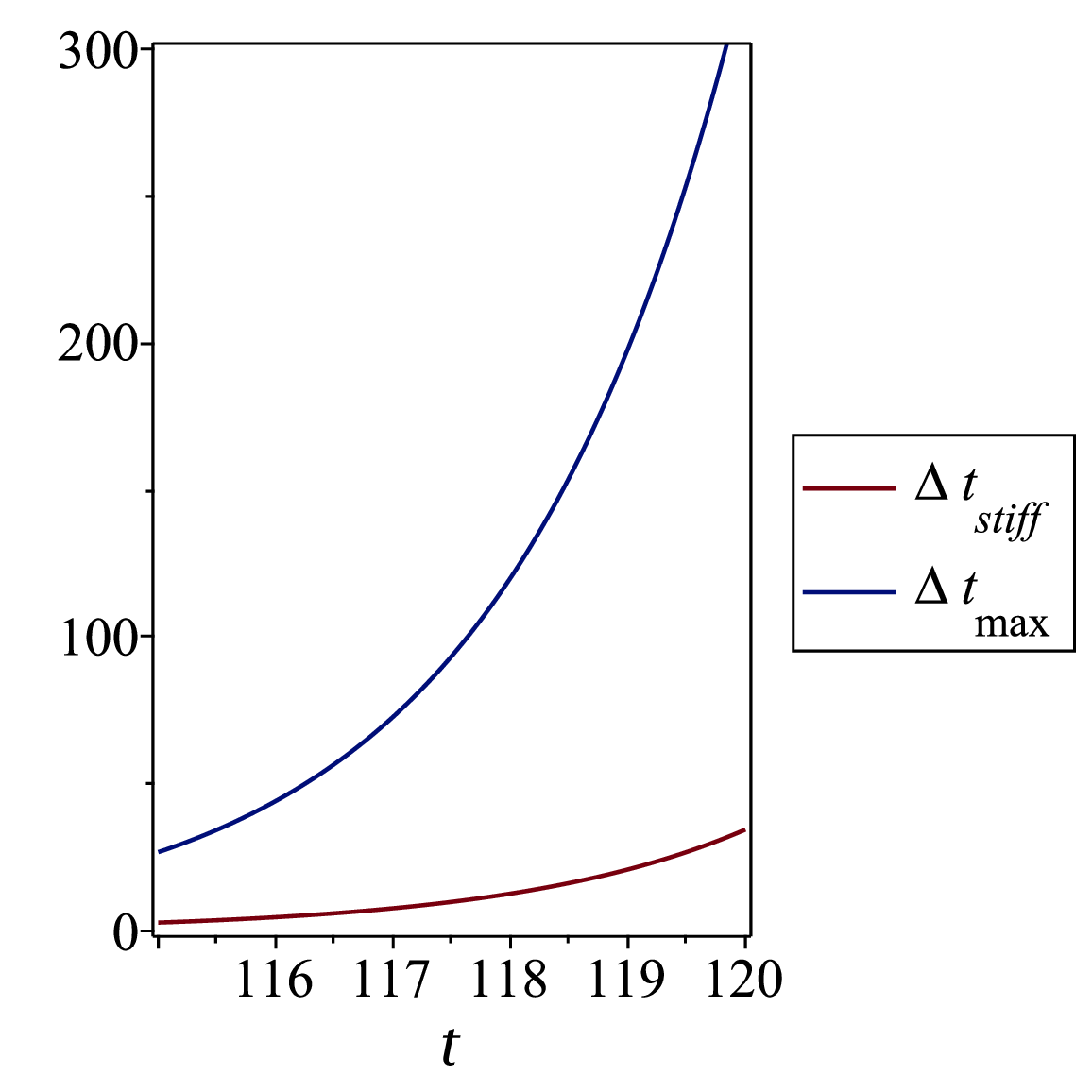}
	\caption{Illustration of stiffness of Eq. (17) for $d=0.01$. Eq. (19) is solved in the two intervals $[105,110]$ and $[115,120]$. The time step for $\epsilon=0.001$ accuracy becomes limited by stiffness sufficiently beyond $t=1/d=100$ where the solution curve $u(t)$ approaches a straight line; thus $(\Delta t)_{stiff}<(\Delta t)_{max}$, or $Q>1$ here.}
\end{figure}

The stiffness of Eq. (17), in the domain beyond $t=1/d$, is demonstrated for $d=0.01$ in Figure 3. We employ $u(t;0)=1$ and $\delta u_0=\epsilon$. The time step for $\epsilon=0.001$ accuracy becomes limited by stiffness ($(\Delta t)_{stiff}<(\Delta t)_{max}$) beyond $t=1/d$ where the solution curve $u(t;0)$ approaches a straight line. 

Explicit solvers are unsuitable when $d<<1$, that is for strong stiffness. Time-spectral ODE solvers, like the GWRM \cite{Scheffel:GWRM1}, are however very efficient; for $d=10^{-4}$, 69 time intervals provide an accuracy of $10^{-4}$ with an efficiency comparable to optimized, implicit Matlab solvers.

\subsection{Strong stiffness - autocatalytical chemical reactions}

\noindent When performing a kinetic analysis of autocatalytic reactions, the following system of equations emerge \cite{Robertson}:

\begin{equation}
	\begin{split}
	\frac{dx}{dt} = -ax+byz  \\
	\frac{dy}{dt} = ax-byz-cy^2  \\
	\frac{dz}{dt} = cy^2 
	\end{split}
\end{equation}

\noindent with $a=0.04$, $b=10^4$, and $c=3\cdot10^7$. The initial conditions are $(x,y,z)=(1,0,0)$. Clearly, the time scales of the reactions are very different. This generally causes no problems for short intervals, but for very long intervals, like $[10^{-6},10^6]$ explicit methods are likely to fail because of the stiffness caused by the large difference between the (reaction rate) constants. As can be seen from the time-adaptive GWRM solution in Figure 4 (note the logarithmic time axis), $y(t)$ features a strong transient initially , whereafter smooth solutions hold for all three variables $x(t),y(t)$ and $z(t)$.

\begin{figure}[h!]
	\centering
	\includegraphics[width=5in]{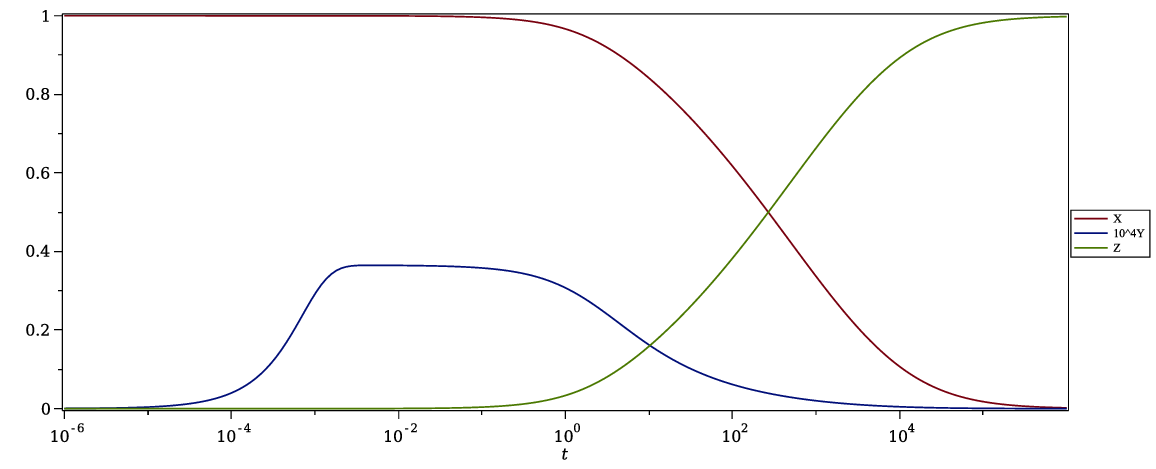}
	\caption{Solution of the Robertson equations (20), employing accuracy $\epsilon=0.001$. Note that $y(t)$ has been magnified a factor $10^4$. The solution was obtained by a GWRM solver.}
\end{figure}

In order to formally study stiffness, we employ the method earlier described. We thus let $\delta u_1 \equiv u_1(t;\epsilon)-u_1(t;0) =  x(t;\epsilon)-x(t;0) \equiv \delta x$ and proceed similarly for the $y$- and $z$- components. Subtraction of the unperturbed and perturbed equations yields, to lowest order in $\delta u_i$;

\begin{equation}
	\begin{split}
	\frac{d \delta x}{dt} = -a\delta x+b(y\delta z + z\delta y)  \\
	\frac{d \delta y}{dt} = a\delta x-b(y\delta z + z\delta y)-2cy\delta y  \\
	\frac{d \delta z}{dt} = 2cy\delta y 
	\end{split}
\end{equation}

The coefficients for $\delta u_i$, including $x,y,z$, are assumed constant in each interval $[t_n,t_{n+1}]$ that is solved for. They do, of course, take on different values in each interval. We now have a system of ODE's, so it is not apparent how the Lyapunov exponents can be determined for deciding stiffness. Since the system is linear, the following theorem may be applied. Assume the system of ODE's

\begin{equation}
	\frac{d\bf u}{dt} = A\bf u
\end{equation}

\noindent where $\textbf u \in R^n$ and $\bf{A}$ is a constant, quadratic, diagonalizable matrix  with distinct eigenvalues ${\gamma}_k \in C$ and eigenvectors $\textbf{c}_k \in C^n$. The solution to Eq. (22) is

\begin{equation}
	\textbf {u}(t) = \sum_{k=1}^na_k e^{\gamma_k} \textbf{c}_k
\end{equation}

\noindent where $a_k$ are arbitrary constants. Comparing with the system of Eqs. (21) we have $\textbf{u} = \delta \textbf{u}$ and $A=J$, where $J$ is the Jacobian of the system. Since we are mainly interested in local Lyapunov exponents for the intervals $[t_n,t_{n+1}]$, we designate the eigenvalues with $\gamma_k$, rather than as $\lambda_k$. These are obtained from the characteristic equation 

\begin{equation}
	|J-\gamma I| = 0
\end{equation}

\noindent where $I$ is the identity matrix. In this case

\begin{equation*}
J = 
\begin{pmatrix}
-a & bz & by \\
a & -bz-2cy & -by \\
0 & 2cy & 0
\end{pmatrix}
\end{equation*}

We may solve Eq. (24) for the eigenvalues numerically. For the initial values, the Robertson equations are not stiff; we obtain $\gamma_1=0$, $\gamma_2=0$ and $\gamma_3=-0.04$ from Eq. (24). But already at $y=10^{-6}$ (maintaining $x=1$ and $z=0$) the system becomes stiff ($\kappa_i$ are large); $\gamma_1=0$, $\gamma_2=-0.05$ and $\gamma_3=-60$. Extreme stiffness is obtained for $t>10^{-4}$, as $y>0.4\cdot10^{-4}$ and $\gamma_3 < -2400$. At this point, $(\Delta t)_{stiff} \approx 2 \sqrt{2}/|\gamma_3| < 10^{-3}$, decreasing further for longer times, and the time step is completely determined by the stiffness of the solution.

An RK4 algorithm with adaptive step size has been employed to solve the problem (20). After over 61 000 time steps the algorithm stagnates near $t=96$. Similar behaviour was encountered for other parameters. It is clear that explicit time-stepping methods are not suited for stiff problems of this type.

Thus an implicit standard trapezoid algorithm with adaptive step size was applied. For an accuracy $\epsilon=0.001$, a solution was obtained, employing 118 steps for the whole interval, with initial step length 0.1. Employing a similar step size adaption algorithm for the same accuracy, the GWRM time-spectral method found the solution about 80 times faster, using 49 time intervals.

\section{Transforming stiff and chaotic differential equations}

\noindent Local Lyapunov stability analysis shows that properties like stiffness and chaoticity depend on the form of the differential equations. Thus, we may ask if it is possible to \textit{transform} the differential equations to produce asymptotically stable solutions. These could then be efficiently and accurately integrated numerically, whereafter a back-transformation is made to the sought solutions. We will find that chaotic equations can be successfully transformed, whereas stiff equations cannot. In the following, this assertion will be demonstrated, employing a general formulation. 

For simplicity we wish to find asymptotically stable solutions to the following single ODE in a time interval $[0,T]$:

\begin{equation}
	\frac{du}{dt} = f(t,u), \quad u(0)=u_0
\end{equation}

\noindent where $f$ is an arbitrary function of the independent variable $t$ and the variable $u(t)$, time-normalized so that $T$ is of order unity. We assume Eq. (25) to be either \textit{stiff} (with neighbouring solutions $\propto exp(\kappa_{0}t)$, with $\kappa_0 < 0$ and $|\kappa_0| >> 1$) or \textit{chaotic} (neighbouring solutions $\delta u(t) \propto exp(\kappa_{0}t)$, with $\kappa_0 > 0$). Strictly, chaotic ODEs feature at least three degrees of freedom, and are thus represented by at least three ODEs [Sprott], but, as will become apparent, the same line of reasoning holds. 

In a local Lyapunov stability analysis we can assume that $\kappa$ is approximately constant in each solution interval $[t_n, t_{n+1}]$. For an asymptotically stable solution, it should hold that $\kappa < 0$ and that $|\kappa t| >> 1$ does not hold.

\subsection{General transformation}

\noindent A general transformation for solving Eq.(25) can be written

\begin{equation}
	u(t) \equiv g(t,z(t))
\end{equation}

\noindent where $g$ may be chosen arbitrarily. In the following, the main focus will be on stiffness. Thus we may, using Eq. (12), assume that $Q_{u}(t)>>1$ at least in some solution intervals, that is for a number of time steps. The goal is now to find an ODE for $z(t)$ for which $Q_{z}(t) \leq 1$ everywhere, guaranteeing asymptotic stability of the solution. From Eqs. (25) and (26) we find

\begin{equation}
	\frac{dz}{dt} = - \frac{\partial g}{\partial t}/\frac{\partial g}{\partial z} + {f(t,g(t,z))}/\frac{\partial g}{\partial z}  \quad .
\end{equation}

\noindent In order to determine

\begin{equation}
	\frac{d\delta z}{dt} = \frac{d(z+\delta z)}{dt} - \frac{dz}{dt}
\end{equation}

\noindent we use the following linearizations

\begin{equation}
	\frac{\partial g}{\partial t} (t, z+ \delta z) = \frac{\partial g}{\partial t} (t, z) + \frac{\partial^2 g}{\partial z \partial t} (t, z)\delta z
\end{equation}

\begin{equation}
	\frac{\partial g}{\partial z} (t, z+ \delta z) = \frac{\partial g}{\partial z} (t, z) + \frac{\partial^2 g}{\partial z^2} (t, z)\delta z
\end{equation}

\begin{equation}
	f(t,g(t,z+\delta z)) = f(t,g(t,z)) + \frac{\partial f}{\partial g} \frac{\partial g}{\partial z} \delta z \quad .
\end{equation}

\noindent Expanding $(\frac{\partial g}{\partial z}(t,z+\delta z))^{-1}$ there results, to first order,

\begin{equation}
	\frac{d\delta z}{dt} = [(\frac{\partial g}{\partial t}-f) \frac{\nicefrac{\partial^2 g}{\partial z^2}}{(\nicefrac{\partial g}{\partial z})^2} - \frac{\nicefrac{\partial^2 g}{\partial t \partial z}}{\nicefrac {\partial  g}{\partial z}} + \frac{\partial f}{\partial g}]\delta z \equiv \kappa_g \delta z  \quad .
\end{equation}

\noindent with the approximate dependence, for $\kappa_g$ $\approx$ constant within a time interval,
\begin{equation}
	\delta z (t) \sim e^{\kappa_g t}
\end{equation}

\noindent We may also write

\begin{equation}
	\frac{d \delta u}{dt} = \frac{\partial f}{\partial u} \delta u \equiv \kappa_f \delta u
\end{equation}

\noindent Now, since

\begin{equation}
	\frac{\partial g}{\partial t} + \frac{\partial g}{\partial z} \frac{dz}{dt} = f,
\end{equation}

\noindent we find from Eq. (32), for $\partial^2g / \partial z^2 \neq 0 $, using $\partial f / \partial g = \partial f / \partial u = \kappa_f$,

\begin{equation}
	\frac{dz}{dt} = (\kappa_f - \kappa_g) \frac {\nicefrac{\partial g}{\partial z}}{\nicefrac{\partial^2 g}{\partial z^2}} - \frac{\nicefrac{\partial^2 g}{\partial t \partial z}}{\nicefrac {\partial^2  g}{\partial z^2}} \quad .
\end{equation}
	
\noindent For an asymptotically stable solution $z(t)$, it holds that $\kappa _g < 0$ and $| \kappa _g \lesssim 1|$. Assuming that Eq. (25) is \textit{stiff}, $\kappa _f < 0$ and $| \kappa _f >> 1|$. Thus $\kappa_g-\kappa_f >>1$ and the ODE in Eq. (36) may be stiff, or strongly chaotic, depending on the choice of $g(t,z(t))$.

\subsection{General transformation, first order}

\noindent The case $\partial^2g / \partial z^2 = 0 $ implies that $g$ is linear in $z$;
\begin{equation}
	g(t,z(t))=A(t)z(t)+B(t) \enskip.
\end{equation}

\noindent Then we, from Eqs. (32) and (34), obtain

\begin{equation}
	\frac{-dA/dt}{A}+\frac{\partial f}{\partial g} = \frac{-dA/dt}{A}+ \kappa_f = \kappa_g ,
\end{equation}

\noindent with formal solution
\begin{equation}
	A(t) = A(0)e^{(\kappa_f - \kappa g) t} .
\end{equation}

\noindent Since $-(\kappa_f-\kappa_g) \approx -\kappa_f >> 1$, $A(t)$ features, as indicated by Eq. (34), a strong local curvature similar to that characterizing the stiffness of Eq. (25). Consequently the ODE for $z(t) = (u(t)-B(t))/A(t)$ will be correspondingly difficult to resolve numerically. The result of the transformation (37) is simply a transition from a stiff ODE for $u(t)$, with typical time step $(\Delta t)^{u}_{stiff}$, to an asymptotically stable ODE for $z(t)$, with typical time step $(\Delta t)^{z}_{max}$, for which it however holds that
\begin{equation}
	\mathcal{O} ((\Delta t)^{z}_{max}) = \mathcal{O} ((\Delta t)^{u}_{stiff})
\end{equation}

This implies that, using the linear transformation (37), neither the linear nor the nonlinear stiff ODEs of Sections 3 and 4 can be transformed to ODEs that can be more efficiently solved for explicit finite difference schemes. Asymptotic stability is attained at the cost of degraded numerical resolution of the sought solution. For the stiff linear ODE in Section 3 we have $\partial f / \partial u = \kappa_f = -a$ so that $A(t) \approx A(0)e^{-at}$. For the stiff nonlinear ODE in Section 4, modelling explosive combustion, the time-normalized ($t=(2/d) \tilde{t}$) relation $\partial f/\partial u = (2/d)(2u-3u^2) \approx -2/d$ holds in the stiff interval $\tilde{t} > 1/2$ where $u \approx 1$ so that $A(\tilde{t} ) \approx A(0)e^{-(2/d)\tilde{t}}$. In both cases $A(t)$, and thus $z(t)$, require high temporal resolution for given accuracy.

\subsection{General transformation, second order}

\noindent Solution of $z(t)$ from Eq. (36) requires that $g(t,z(t))$ is expanded to second order. We obtain

\begin{equation}
	\frac{\partial g}{\partial z} \approx A + B(t-t_0) + C(z-z_0)
\end{equation} 

\noindent where
\begin{equation}
	A \equiv \frac{\partial g}{\partial z}(t_0,z_0), B \equiv \frac{\partial^2 g}{\partial t \partial z}(t_0,z_0), C \equiv \frac{\partial^2 g}{\partial z^2}(t_0,z_0) .
\end{equation} 

\noindent Considering individual time intervals $[t_n,t_{n+1}]$, we can let $t_0 = t_n = 0$ without losing generality. The notation $z_0 = z(t_0) = z(0)$ is also used. Inserting Eq. (41) into Eq. (36), the following solution is obtained, exact to second order;
\begin{equation}
	z(t) = \frac{A}{C}(e^{(\kappa_f-\kappa_g)t}-1) - \frac{B}{C} t + z_0 .
\end{equation} 

\noindent Recognizing that $-(\kappa_f-\kappa_g) >> 1$, the temporal resolution required is again approximately the same as that required due to the stiffness of the problem. 

In order to also study the effect of explicit $t$-dependence in $g(t,z(t))$, the following transformation is useful:  
\begin{equation}
	g(t,z(t)) = g_0 e^{\kappa _z t} z^\alpha
\end{equation} 

\noindent where $g_0$, $\kappa_z$ are constants and $\alpha >1$. Linearizing in $t$, using Eqs. (41)-(43), there obtains
\begin{equation}
	z(t) \approx z(0) e^{((\kappa_f - \kappa_g - \kappa_z)/(\alpha-1))t}  ,
\end{equation}

\noindent Choosing $\kappa_z$ so that $\kappa_f - \kappa_g - \kappa_z$ becomes small, requires that $\kappa_z \approx \kappa_f$, which, according to Eq. (44), violates the assumption that $u(t) = g(t,z(t))$ should be a smooth function of $t$.

That explicit time dependence in $g(t,z(t))$ cannot compensate for the effect of the term $\kappa_f-\kappa_g$ in $z(t)$ can also be seen by setting $dz/dt \approx 0$ in Eq. (36). It follows that  
\begin{equation}
	(\kappa_f - \kappa_g) \frac {\partial g}{\partial z} \approx \frac {\partial^2 g}{\partial t \partial z}
\end{equation}

\noindent The solution to this equation is obtained by integration in $t$ and $z$;

\begin{equation}
	g(t,z(t)) = e^{(\kappa_f-\kappa_g)t} ( \int_{0}^t e^{-(\kappa_f-\kappa_g)t'}T(t')dt' + Z(z))
\end{equation}

\noindent where $T(t)$ and $Z(z)$ are arbitrary functions. Since we may assume $Z(z) \neq 0$, the last term will be non-smooth and again the problem requires a numerical time resolution amounting to that needed for resolving stiffness. 

If the problem (25) on the other hand were \textit{chaotic}, so that $\kappa_f>0$, we could maintain the asymptotically stable conditions $\kappa _g < 0$ and $| \kappa _g \lesssim 1|$ for the numerical solution $z(t)$ by a suitable choice of $\kappa_z$. It should be noted that the function $g(t,z(t))$ must be assumed given in order to solve for $z(t)$.

\section{Removing numerical chaoticity}

\noindent The question thus arises how the numerical errors caused by chaotic differential equations can be mitigated locally, using variable transformations. The idea is to, via a suitable transformation, temporarily remove the exponential divergence between neighbouring solutions in parameter space for the new variable. The sensitivity to initial conditions at each new time step is thus alleviated. Due to the nonlinearity of the equations, the net result is substantially higher accuracy than what can be attained by standard time-stepping solvers like explicit RK4. We will now illustrate how accuracy can be substantially improved for a well known set of chaotic ODEs.

\subsection{The Lorenz 1984 chaotic equations}

\noindent The Lorenz 1984 model \cite{Lorenz:1984atmosphere} is a set of three non-linear ordinary differential equations that features chaotic, rather than stiff, behaviour related to meteorological systems. It is one of the simplest models of Hadley circulation. Thus, whereas the Lorenz 1984 model is not an accurate numerical weather prediction (NWP) model, it is of interest for numerical analysis of simple, yet non-linear chaotic behaviour that may occur in NWP. The equations are
\begin{align}
\frac{dx}{dt} &= -y^2-z^2-ax+aF, \label{eq:L1} \\
\frac{dy}{dt} &= xy-bxz-y+G, \nonumber \\
\frac{dz}{dt} &= bxy+xz-z \nonumber
\end{align}
Further details of the model are found in \cite{Lorenz:1984atmosphere}. We note here that the variables $x(t)$, $y(t)$, and $z(t)$ represent certain meteorological systems such as wind currents and large-scale eddies. The coefficients $a$, $b$, $F$, and $G$ are chosen within certain limits to act as damping, coupling, and amplification of the physical processes of the system. In Figure 5, a solution for $a=0.25$, $b=4.0$, $F=8.0$ and $G=1.0$ is shown. The initial conditions are $(x,y,z)=(0.96,-1.1,0.5)$.

\begin{figure}[h!]
	\centering
	\includegraphics[width=5in]{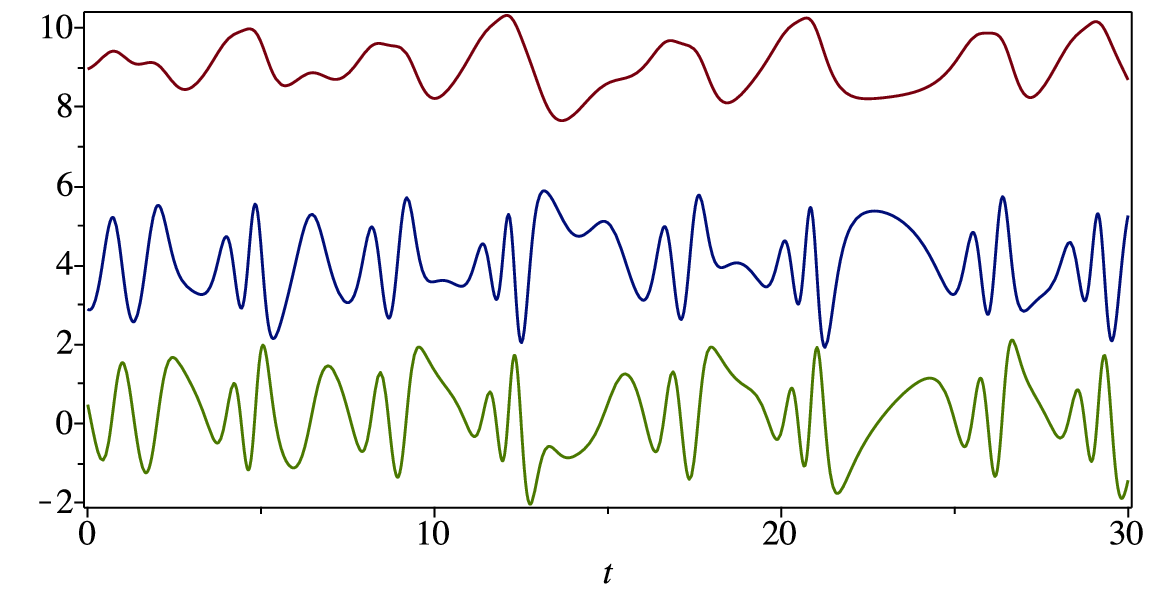}
	\caption{Solution of the Lorenz equations (48) with parameters $a=0.25$, $b=4.0$, $F=8.0$ and $G=1.0$. The initial conditions are $(x,y,z)=(0.96,-1.1,0.5)$. From top to bottom; $x(t)+8$, $y(t)+4$, $z(t)$. A high accuracy, time-spectral GWRM code was used.}
\end{figure}

This set of equations is difficult to solve for any explicit time-stepping method, since chaoticity causes any error in the solutions at time $t$ to be strongly amplified in the interval $[t_n,t_{n+1}]$. To see this, we analyze the eigenvalues as in the previous section, thus obtaining the system

\begin{equation}
	\begin{split}
	\frac{d \delta x}{dt} = -2y\delta y -2z\delta z - a\delta x  \\
	\frac{d \delta y}{dt} = x\delta y + y\delta x - bx\delta z - bz\delta x - \delta y  \\
	\frac{d \delta z}{dt} = bx\delta y + by\delta x + x \delta z + z\delta x - \delta z 
	\end{split}
\end{equation}

\noindent The corresponding Jacobian becomes

\begin{equation}
J = 
\begin{pmatrix}
-a & -2y & -2z \\
y-bz & x-1 & -bx \\
by+z & bx & x-1
\end{pmatrix}
\end{equation}

\noindent For the initial values, the Lorenz equations are chaotic. By solving the characteristic equation $|J-\gamma I| = 0$ for $t=0$ there obtains $\gamma_1=1.9$, $\gamma_2=-1.1+4.5i$ and $\gamma_3=-1.1-4.5i$. It is seen from Eq. (49) that the parameters $F$ and $G$ do not affect the character of the ODE's. The system is found to be robust in the sense that chaoticity is always manifested by at least one eigenvalue being positive for any value of $b$. For $a \geq 3.1$, however, the system will initially be non-chaotic and may also locally be non-chaotic, for example for $t \approx 13-14$. The equations are not stiff for any values of $a$ and $b$. Stiffness and chaoticity are indeed different properties.

It is of interest to determine the local eigenvalues of the characteristic equation for the full solution interval $t \in [0,30]$. In Figure 6, these have been computed at 400 equidistant times. It is seen that at least one local Lyapunov exponent is positive, indicating chaoticity, in the entire interval except for in small intervals near $t=14$.

\begin{figure}[h!]
	\centering
	\includegraphics[width=5in]{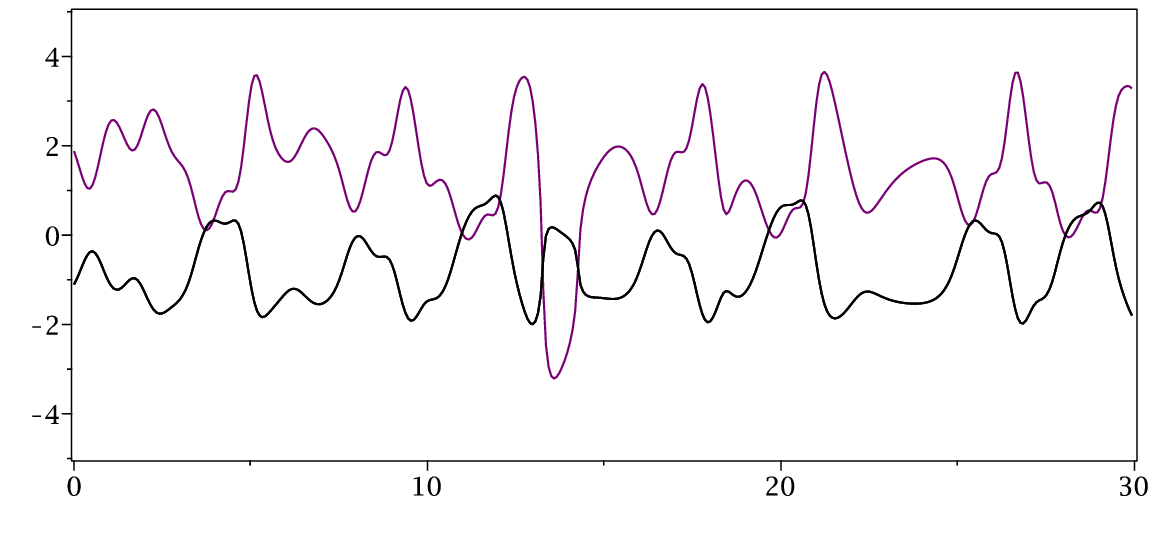}
	\caption{Eigenvalues $\gamma$ of the characteristic equation $|J-\gamma I| = 0$ for the Lorenz local Lyapunov exponent system (49) in the time interval $[0,30]$. The top curve represents the real root, whereas the bottom curve images the real part of the remaining two, complex conjugated roots.}
\end{figure}

\subsection{Transforming the Lorenz equations to eliminate numerical chaoticity}
\noindent Recognizing the exponential dependence of chaoticity we choose the following, simple transformations for the Lorenz equations, where $A_i(t)$ are explicit functions of $t$, and $\epsilon_i$ as well as $\mu_i$ are constants:

\begin{equation}
	\begin{split}
	x(t)=A_1(t)z_1(t)=\epsilon_1 e^{\mu_1t} z_1(t)  \\
	y(t)=A_2(t)z_2(t)=\epsilon_2 e^{\mu_2t} z_2(t)  \\
	z(t)=A_3(t)z_3(t)=\epsilon_3 e^{\mu_3t} z_3(t) 
	\end{split}
\end{equation}

Solving for the $z$ variables, Eq. (48) becomes

\begin{equation}
	\begin{split}
	\frac{dz_1}{dt}=-\mu_1 z_1 - \frac{\epsilon_2^2}{\epsilon_1} e^{(2\mu_2-\mu_1)t} z_2^2  - \frac{\epsilon_3^2}{\epsilon_1} e^{(2\mu_3-\mu_1)t} z_3^2-az_1+\frac{a}{\epsilon_1}e^{-\mu_1 t}F\\
	\frac{dz_2}{dt}=-\mu_2 z_2 + \epsilon_1 e^{\mu_1 t}z_1 z_2    - b\frac{\epsilon_1 \epsilon_3}{\epsilon_2} e^{(\mu_1 -\mu_2+\mu_3)t} z_1 z_3-z_2+\frac{G}{\epsilon_2}e^{-\mu_2 t}\\
	\frac{dz_3}{dt}=-\mu_3 z_3 + b\frac{\epsilon_1 \epsilon_2}{\epsilon_3} e^{(\mu_1+\mu_2 -\mu_3)t} z_1 z_2 +\epsilon_1 e^{\mu_1}z_1 z_3 -z_3   \\
	\end{split}
\end{equation}

\noindent This system will be solved with a standard RK4 algorithm, using $N$ equidistant time steps for $t \in [0,T]$. In order to avoid excessively large numerical exponential factors, the time domain is separated into $K$ intervals with $N/K$ time steps in each. In all cases to follow $N=600$. At the end of each interval, new initial conditions for $z_i(t)$ are generated using Eq. (51). For simplicity, we here set $\epsilon_1 = \epsilon_2 = \epsilon_3 = 1$. 

\subsubsection{Method 1}

\noindent The problem is now to optimize the choice of $\mu_i$ in order to obtain asymptotically stable, rather than chaotic, solutions $z_i(t)$. Given the knowledge of the local Lyapunov exponents displayed in Figure 6, one may be tempted to simply set all $\mu_i=\mu_0$, where $\mu_0$ is a positive constant, for example in the interval $[2,3]$. This approach may not be effective, since optimal choices of $\mu_i$ are expected to be interdependent. A number of RK4 code runs were run manually, individually varying $\mu_i$. To compute the local Lyapunov exponents for the transformed equations (52), we determine, as before, the system 
\begin{equation}
	\begin{split}
	\frac{d \delta z_1}{dt} = -\mu_1 \delta z_1 - 2z_2\delta z_2 - 2z_3\delta z_3 - a\delta z_1  \\
	\frac{d \delta z_2}{dt} = -\mu_2 \delta z_2 + z_1\delta z_2 + z_2\delta z_1 - bz_1\delta z_3 - bz_3\delta z_1 - \delta z_2  \\
	\frac{d \delta z_3}{dt} = -\mu_3 \delta z_3 + bz_1\delta z_2 + bz_2\delta z_1 + z_1 \delta z_3 + z_3\delta z_1 - \delta z_3 
	\end{split}
\end{equation} 

We make the approximation $t=0$ in Eq. (52), corresponding to short interval lengths $t_{n+1}-t_n$. The eigenvalues (local Lyapunov exponents) $\gamma$ are thus obtained from the determinant $|J^*-\gamma I| = 0$ where

\begin{equation}
J^* = 
\begin{pmatrix}
-a - \mu_1 & -2z_2 & -2z_3 \\
z_2-bz_3 & z_1-1 -\mu_2 & -bz_1 \\
bz_2+z_3 & bz_1 & z_1-1 - \mu_3
\end{pmatrix}
\end{equation}

In Figure 7 the eigenvalues $\gamma$ for the choice $\mu_1=2.592$, $\mu_2=1.944$ and $\mu_3=1.539$ (constant throughout the computations) are shown. It is seen that the local Lyapunov exponents have been reduced substantially, so that the system (52), solved with retained explicit time dependence for accuracy, becomes non-chaotic over a large fraction of the time domain $[0,30]$. The ratio $N/K=10$ (10 time steps in each time interval) was chosen. We term this approach \textit{Method 1}. Three additional methods are discussed in the following, whereafter all methods will be compared in Figure 8.

\begin{figure}[h!]
	\centering
	\includegraphics[width=5in]{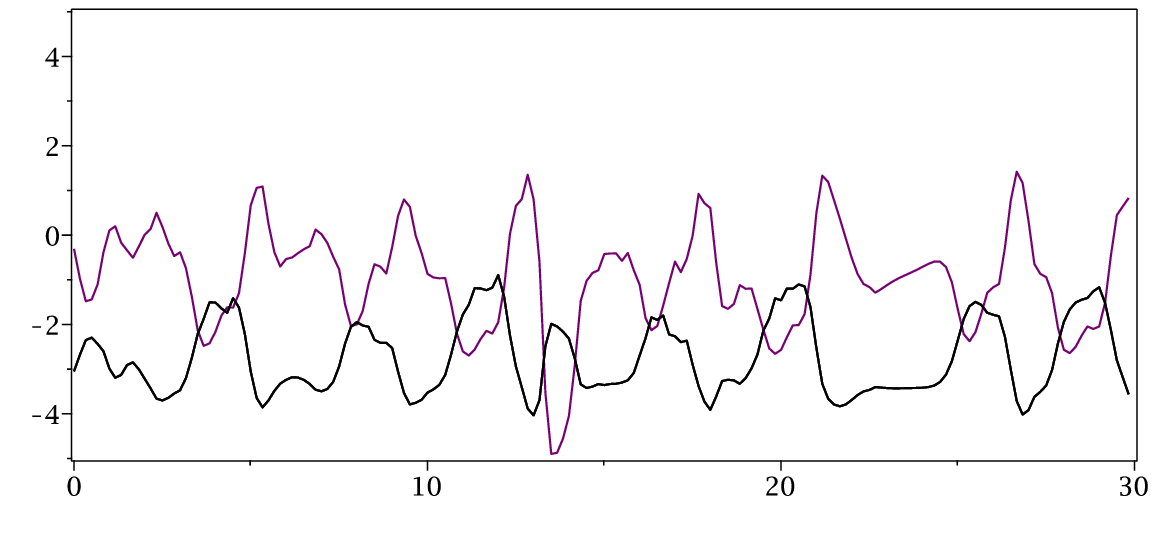}
	\caption{Eigenvalues $\gamma$ of the characteristic equation $|J^*-\gamma I| = 0$ for the Lorenz local Lyapunov exponent system (53) with $\epsilon_1 = \epsilon_2 = \epsilon_3 = 1$, using $\mu_1=2.592$, $\mu_2=1.944$ and $\mu_3=1.539$  in the time interval $[0,30]$. The top curve represents the real root, whereas the bottom curve images the real part of the remaining two, complex conjugated roots. A large fraction of the temporal domain is non-chaotic.}
\end{figure}

\subsubsection{Method 2}

\noindent It may be noted that, in the determinant $|J^*-\gamma I| = 0$, the eigenvalues and the transformation parameters will occur as terms $\mu_i+\gamma$ only.  This shows that i) the chosen transformation (51) features a high degree of simplicity and interpretability and ii) that $\mu_i$ typically should be chosen sufficiently positive to reduce the real parts of the eigenvalues $\gamma$ towards negative values, corresponding to asymptotically stable solutions. A complication is that the interrelation between the constants $\mu_i$ and the eigenvalues is unclear, since for each choice of $(\mu_1,\mu_2,\mu_3)$ there generally exists three eigenvalues $\gamma$. 

An option would be to simply set all eigenvalues equal (multiple roots), to a suitable negative value for asymptotic stability and solve for $(\mu_1,\mu_2,\mu_3)$. The corresponding, multiple root solutions would then be of the form $\delta z_i = (c_{2i}t^2+c_{1i}t+c_{0i})e^{\gamma t}$, where the $c$'s are constants. Assuming short time intervals, $\delta z_i \approx c_{0i}e^{\gamma t}$ where the time dependence is governed by the exponential (Lyapunov-like) term alone. This approach is, however, only moderately successful with respect to minimizing the error for longer times ($t>10$). This is due to the rapid time variation of the eigenvalues, as can be seen in Figs. 6 and 7; knowledge of the LLEs is not sufficiently helpful for providing asymptotic stability in the subsequent time interval.

Method 2 is an algorithm, that instead compensates for local chaoticity. The maximum real part of the three eigenvalues $\gamma_{max}(t_{n})$ of $|J^*-\gamma I| = 0$ is computed in each time interval $[t_n,t_{n+1}]$, and the values $\mu_i = q \gamma_{max}(t_{n})$, where $q$ is a constant, are used for the subsequent interval. For the case shown here, the initial values were $\mu_1=2.0, \mu_2=2.0, \mu_3=2.0$ and $q=1.0$ was used for $N/K=10$. 

\subsubsection{Methods 3 and 4}

\noindent In Method 3, $\mu_i$ values for subsequent time intervals are computed from accumulated averages $<\gamma_{max}>(t_{n})$ of the maximum eigenvalues $\gamma_{max}(t_{n})$. High accuracy is obtained for the choice $\mu_1(t_{n})=-1.5<\gamma_{max}>(t_{n}), \mu_2(t_{n})=-0.66<\gamma_{max}>(t_{n}), \mu_3(t_{n})=-0.5<\gamma_{max}>(t_{n})$, employing $\mu_1(t_{1})=2.16, \mu_2(t_{1})=1.62, \mu_3(t_{1})=1.28$. 

High accuracy also results from Method 4, where averages $<\gamma_{max}>(t_{n})$ are taken from the last two previous time intervals only. Both methods employ 40 time steps in each time interval, that is $N/K=40$.

\subsubsection{Method comparisons}

\noindent In Figure 8 the absolute errors for $x(t)$ of Eq. (48), computed by the transformation Methods 1-4, are compared with that of standard RK4 solution. It is seen that radical improvements in accuracy is obtained using the transformations. Methods 3 and 4 are about 100 times more accurate in the entire interval than standard RK4, and display only a weak deterioration with time. In Figure 9 the exact $x(t)$ solution, with same parameters as in Figure 5, is compared to that obtained from standard RK4 and Method 3, for the numerical parameters given earlier.  

It should be noted from Figure 8 that although RK4 accuracy is substantially improved by the four transformation methods discussed, asymptotic stability (error tending to zero) is not achieved since the Lorenz 1984 equations are inherently chaotic; any error in the solution will be amplified for higher times even for infinite numerical accuracy. The effect of the transformations is to render the $z(t)$ solutions asymptotically stable, thereby nonlinearly mitigating errors in the numerical, explicit time-stepping algorithm. 

\begin{figure}[h!]
	\centering
	\includegraphics[width=5in]{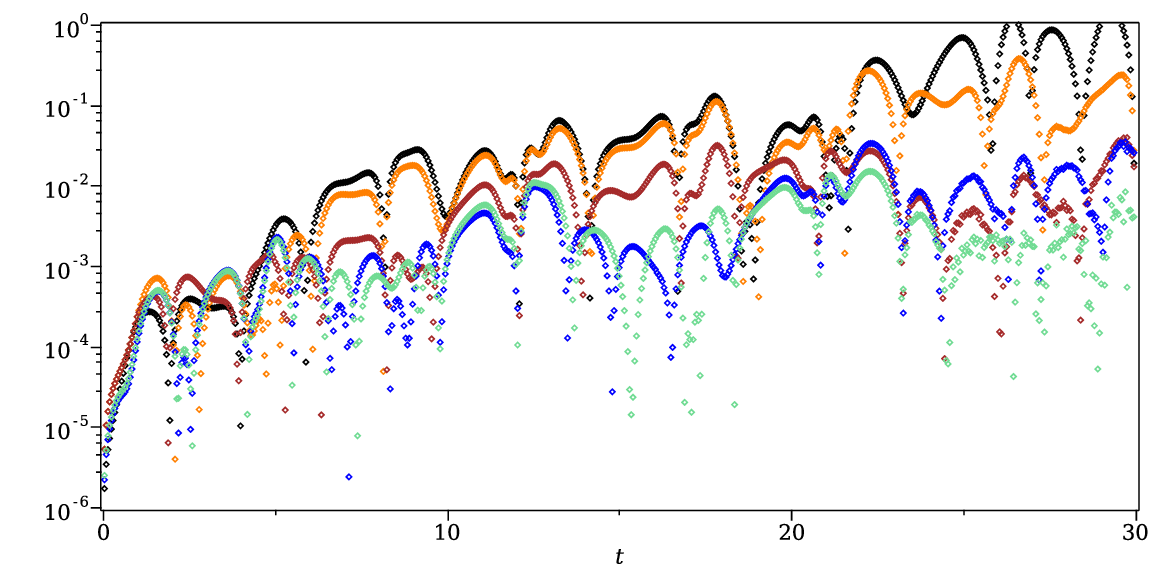}
	\caption{Absolute errors for RK4 solution of $x(t)$ from the Lorenz equations (48), with same parameters as in Figure 5, using standard RK4 (top, black points) and for 4 transformation methods, as described in the text: Method 1 (orange), Method 2 (brown), Method 3 (blue), Method 4 (green).}
\end{figure}

\begin{figure}[h!]
	\centering
	\includegraphics[width=5in]{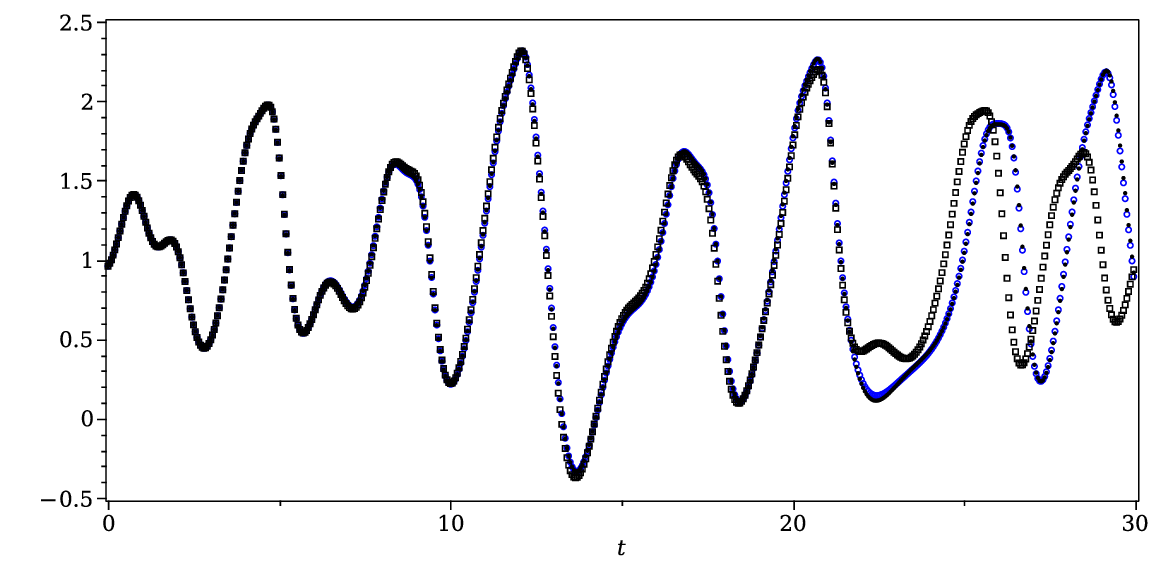}
	\caption{$x(t)$ obtained from the Lorenz equations (48), with Figure 5 parameters: exact solution (points), standard RK4 (squares), and Method 3 (blue).}
\end{figure}

\subsubsection{Time step extension}
\noindent We now seek to investigate whether our objective has been achieved, namely, that the influence of chaotic behavior has been reduced to the point where the achieved level of accuracy is primarily constrained by numerical approximation thresholds rather than by time step limitations caused by chaotic dynamics. To investigate this, we compare the local $(\Delta t)_{max}(t)$, given by Eq. (7) and the curvature of the exact solution, with the step size $\Delta= 30/600=0.05$, employed in the present computations. 

The curvature is computed from the most non-smooth solution $Z(t)$, rather than from $X(t)$ or $Y(t)$. We use the maximum error $(\epsilon=0.012)$ of Method 3 in the time interval $[0,30]$ for computing $(\Delta t)_{max}(t)$.  It is clear, from Figure 10, that $\Delta$, being fixed in the interval, could not have been longer, since then the requirement $(\Delta t)_{max}(t) \geq \Delta$ would be locally violated. Inserted is also $(\Delta t)_{max}(t)$ as obtained from standard RK4, for which $\epsilon=1.0$. Clearly, the limited accuracy for standard RK4 is determined by chaoticity, rather than by $(\Delta t)_{max}(t)$. Thus we may conclude that the transformations of Methods 3 and 4 remove the time step limitations that are associated with the chaoticity of the Lorenz equation (48).  

This conclusion is further supported by Figure 11, where Method 3 is applied using a larger number of time steps $N=1620$ ($\Delta=0.0185$) with $K=60$, overall obtaining an accuracy $\epsilon=0.001$. It is again seen that longer $\Delta$ would violate the maximum allowable time step for numerical resolution $(\Delta t)_{max}(t)$ at this accuracy. 

\begin{figure}[h!]
	\centering
	\includegraphics[width=5in]{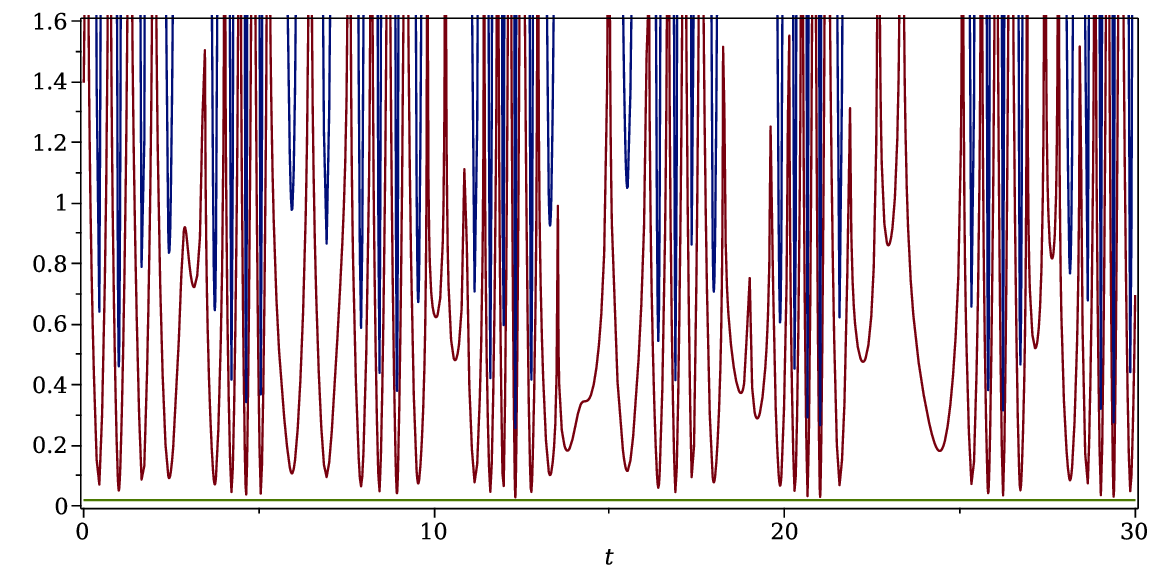}
	\caption{Comparison of $(\Delta t)_{max}(t)$, given by Eq. (7) and obtained from the curvature of the exact solution for the overall accuracy $\epsilon=0.012$ obtained by Method 3, with step size $\Delta=30/600=0.05$. The curvature is computed from $Z(t)$, being the most non-smooth solution. Inserted is also $(\Delta t)_{max}(t)$ (in blue) obtained by standard RK4, where $\epsilon=1.0$. }
\end{figure}

\begin{figure}[h!]
	\centering
	\includegraphics[width=5in]{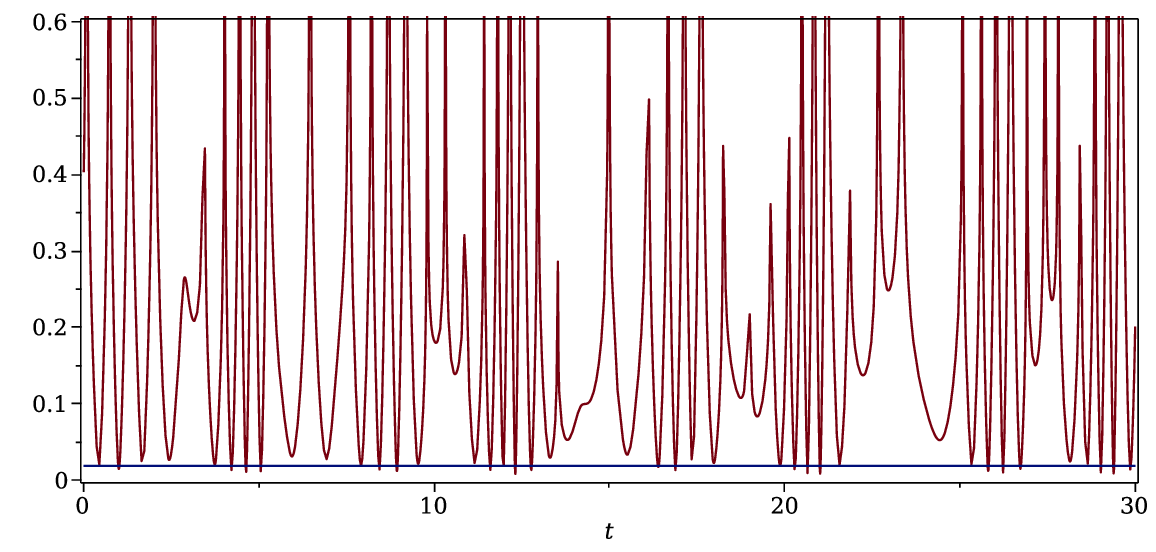}
	\caption{$(\Delta t)_{max}(t)$ as in Figure 10, but for Method 3 only, employing overall accuracy $\epsilon=0.001$. Obtained for $N=1620$ ($\Delta=0.0185$, also displayed in the figure) with $K=60$.}
\end{figure}

\section{Discussion}

\noindent Chaoticity, when appearing in systems of differential equations, is inherent and cannot be transformed away. But the numerical errors introduced by the discreteness of time-stepping methods can be mitigated by transformations to systems that produce asymptotically stable solutions. In these systems, the tendency for neighbouring solutions to converge to each other compensates for the sensitivity to initial conditions that characterize chaotic systems. In Section 6.2 we suggest an analytically simple transformation method along these lines. Based on local Lyapunov exponents, it is applicable to time intervals consisting of a single or a number of time steps. At the end of each interval, a back transformation to the original variables are used, producing new initial conditions. 

The back transformation involves an exponential amplification of the numerical errors but, as we have seen by solving the chaotic Lorenz 1984 ODEs, the net result is substantially higher accuracy than what can be attained by standard time-stepping solvers like explicit fourth order Runge-Kutta. Employing four different methods for determining the transformation parameters $\mu_i$, it was found that the time step length was no longer determined by chaoticity but by discrete, numerical resolution of the exact solution at the given accuracy.

Interesting topics for further research is to find robust methods for determining parameters $\mu_i$ that result in accurate solutions for a variety of chaotic problems. The transformation of Section 6.2 should be generalizable also to PDEs. Thus, problems in areas like numerical weather prediction and turbulence, which for long have been plagued by time step limitations, could be addressed.

Turning to stiff systems, it is no major surprise that these problems cannot be remedied by transformations. In the formal analysis presented it also becomes clear, by analysis of the ODEs governing the local Lyapunov exponents, that the original difficulty to numerically resolve neighbouring, fast declining solutions, is transferred to the numerical resolution of the transformed functions. The required analytical tools $(\Delta t)_{stiff}$ and $(\Delta t)_{max}$ were introduced in Section 2.

Fixed time step lengths are used in this work in order to isolate the effect of the transformations alone. It is expected that adaptive step length algorithms, in combination with suitable methods for choosing transformation parameters (like $\mu_i$) could help to optimize the effect of the transformation, thus enabling even longer time steps for solution of chaotic ODEs and PDEs with retained accuracy. 

\section{Conclusion}
\noindent Transform methods for improving accuracy of explicit time-stepping ODE solvers are introduced. It is shown that the numerical accuracy of solutions to chaotic ODEs can be improved by at least two orders of magnitude, at retained step length, by transforming to ODEs with asymptotically stable solutions. The chaotic Lorenz 1984 equations serves as an example.

For stiff problems, however, accuracy remains a problem when transforming to problems with asymptotically stable solutions. Both linear and nonlinear examples of stiff systems are studied. The results rest on the concept of local Lyapunov exponents (LLEs), for which we have defined a measure $Q(t)$ for diagnosing asymptotic stability, stiffness and chaos. Applications to stiff linear and nonlinear ODEs are given. 


\label{}





\bibliographystyle{elsarticle-num}
\bibliography{mybib}

\begin{thebibliography}{10}
\expandafter\ifx\csname url\endcsname\relax
  \def\url#1{\texttt{#1}}\fi
\expandafter\ifx\csname urlprefix\endcsname\relax\def\urlprefix{URL }\fi
\expandafter\ifx\csname href\endcsname\relax
  \def\href#1#2{#2} \def\path#1{#1}\fi

\bibitem{Hirsch}
M.~W. Hirsch, S.~Smale, Differential Equations, Dynamical Systems and Linear
  Algebra, Academic, New York, 1974.

\bibitem{Sprott}
J.~C. Sprott, Some simple chaotic jerk functions, Am. J. Phys. 65~(6) (1997)
  537--543.

\bibitem{Cartwright:1}
J.~H. Cartwright, {Nonlinear stiffness, Lyapunov exponents, and attractor
  dimension}, Physics Letters A 264 (1999) 298--302.

\bibitem{Lambert}
J.~Lambert, Numerical Methods for Ordinary Differential Systems, New York:
  Wiley, 1992.

\bibitem{Scheffel:GWRM2}
J.~Scheffel, {A Spectral Method in Time for Initial-Value Problems}, American
  Journal of Computational Mathematics 02 (2012) 173--193.
\newblock \href {https://doi.org/10.4236/ajcm.2012.23023}
  {\path{doi:10.4236/ajcm.2012.23023}}.

\bibitem{Heath}
M.~T. Heath, Scientific Computing - an Introductory Survey, McGraw-Hill, 2005.

\bibitem{Hairer}
E.~Hairer, G.~Wanner, Solving Ordinary Differential Equations II: Stiff and
  Differential - Algebraic Problems, Springer, 2nd ed. 1996.

\bibitem{Ayers}
D.~Ayers, J.~Lau, J.~Amezcua, A.~Carrassi, V.~Ojha,
  \href{https://arxiv.org/abs/2202.04944}{{Supervised machine learning to
  estimate instabilities in chaotic systems: estimation of local Lyapunov
  exponents}}, arXiv, (2022).
\newblock \href {https://doi.org/10.48550/ARXIV.2202.04944}
  {\path{doi:10.48550/ARXIV.2202.04944}}.
\newline\urlprefix\url{https://arxiv.org/abs/2202.04944}

\bibitem{Scheffel:GWRM3}
J.~Scheffel, K.~Lindvall, H.~F. Yik, {A time-spectral approach to numerical
  weather prediction}, Computer Physics Communications 226 (2017) 127--135.
\newblock \href {https://doi.org/10.1016/j.cpc.2018.01.010}
  {\path{doi:10.1016/j.cpc.2018.01.010}}.

\bibitem{Scheffel:GWRM1}
J.~Scheffel, Time-spectral solution of initial-value problems, in: Partial
  Differential Equations: Theory, Analysis and Applications, Nova Science
  Publishers, Inc., 2011, pp. 1--49.

\bibitem{Robertson}
H.~H. Robertson, The solution of a set of reaction rate equations, in: J.~Walsh
  (Ed.), Numerical Analysis: An Introduction, Academic Press, London, 1966, pp.
  178--182.

\bibitem{Lorenz:1984atmosphere}
E.~N. Lorenz, {Irregularity: a fundamental property of the atmosphere}, Tellus
  A 36 (1984) 98--110.
\newblock \href {https://doi.org/10.3402/tellusa.v36i2.11473}
  {\path{doi:10.3402/tellusa.v36i2.11473}}.

\end{thebibliography}


\begin{thebibliography}{0}
\bibitem{1}Reference 1         
\bibitem{2}Reference 2         
\bibitem{3}Reference 3         
\end{thebibliography}







\end{document}